\documentclass{amsart}
\usepackage{amsmath,amssymb,hyperref}
\usepackage{amsrefs}
\usepackage{mathrsfs}
\newtheorem{thm}{Theorem}

\newtheorem{cor}[thm]{Corollary}
\newtheorem{lemma}[thm]{Lemma}
\newtheorem{prop}[thm]{Proposition}

\theoremstyle{definition}
\newtheorem{definition}[thm]{Definition}

\theoremstyle{remark}
\newtheorem{remark}[thm]{Remark}
\newtheorem{example}[thm]{Example}
\def\mathcs{C^{*}}
\newcommand{\cs}{\ensuremath{\mathcs}}

\DeclareMathSymbol{\rtimes}{\mathbin}{AMSb}{"6F}
\newcommand{\ib}{im\-prim\-i\-tiv\-ity bi\-mod\-u\-le}

\newcommand{\sme}{\,\mathord{\mathop{\text{--}}\nolimits_{\relax}}\,}

\def\K{\mathcal{K}}

\DeclareMathOperator{\Ind}{Ind}

\DeclareMathOperator*{\supp}{supp}
\def\set#1{\{\,#1\,\}}
\newcommand\sset[1]{\{#1\}}
\let\tensor=\otimes
\def\restr#1{|_{{#1}}}
\makeatletter
\def\labelenumi{\textnormal{(\@alph\c@enumi)}}
\def\theenumi{\@alph \c@enumi}
\def\labelenumii{\textnormal{(\@roman\c@enumii)}}
\def\theenumii{\@roman \c@enumii}
\newcount\charno
\def\alphapart#1{\charno=96
\advance\charno by#1\char\charno}

\makeatother

\def\<{\langle}
\def\>{\rangle}
\let\ipscriptstyle=\scriptscriptstyle
\def\lipsqueeze{{\mskip -3.0mu}}
\def\ripsqueeze{{\mskip -3.0mu}}
\def\ipcomma{\nobreak\mathrel{,}\nobreak}
\newbox\ipstrutbox
\setbox\ipstrutbox=\hbox{\vrule height8.5pt
width 0pt}
\def\ipstrut{\copy\ipstrutbox}
\def\lip#1<#2,#3>{\mathopen{\relax_{\ipstrut\ipscriptstyle{
#1}}\lipsqueeze
\langle} #2\ipcomma #3 \rangle}
\def\blip#1<#2,#3>{\mathopen{\relax_{\ipstrut
\ipscriptstyle{ #1}}\lipsqueeze\bigl\langle} #2\ipcomma #3 \bigr\rangle}
\def\rip#1<#2,#3>{\langle #2\ipcomma #3
\rangle_{\ripsqueeze\ipstrut\ipscriptstyle{#1}}}
\def\brip#1<#2,#3>{\bigl\langle #2\ipcomma #3
\bigr\rangle_{\ripsqueeze\ipstrut\ipscriptstyle{#1}}}
\def\angsqueeze{\mskip -6mu}
\def\smangsqueeze{\mskip -3.7mu}
\def\trip#1<#2,#3>{\langle\smangsqueeze\langle #2\ipcomma #3
\rangle\smangsqueeze\rangle_{\ripsqueeze\ipstrut\ipscriptstyle{#1}}}
\def\btrip#1<#2,#3>{\bigl\langle\angsqueeze\bigl\langle #2\ipcomma
#3
\bigr\rangle
\angsqueeze\bigr\rangle_{\ripsqueeze\ipstrut\ipscriptstyle{#1}}}
\def\tlip#1<#2,#3>{\mathopen{\relax_{\ipstrut\ipscriptstyle{
#1}}\lipsqueeze \langle\smangsqueeze\langle} #2\ipcomma #3
\rangle\smangsqueeze\rangle}
\def\btlip#1<#2,#3>{\mathopen{\relax_{\ipstrut\ipscriptstyle{
#1}}\lipsqueeze
\bigl\langle\angsqueeze\bigl\langle} #2\ipcomma #3
\bigr\rangle\angsqueeze\bigr\rangle}

\def\ip(#1|#2){(#1\mid #2)}
\def\bip(#1|#2){\bigl(#1 \mid #2\bigr)}
\def\Bip(#1|#2){\Bigl( #1 \bigm| #2 \Bigr)}
\expandafter\ifx\csname BibSpec\endcsname\relax\else
\BibSpec{collection.article}{%
    +{}  {\PrintAuthors}                {author}
    +{,} { \textit}                     {title}
    +{.} { }                            {part}
    +{:} { \textit}                     {subtitle}
    +{,} { \PrintContributions}         {contribution}
    +{,} { \PrintConference}            {conference}
    +{}  {\PrintBook}                   {book}
    +{,} { }                            {booktitle}
    +{,} { }                            {series}
    +{,} { \voltext}                    {volume}
    +{,} { }                            {publisher}
    +{,} { }                            {organization}
    +{,} { }                            {address}
    +{,} { \PrintDateB}                 {date}
    +{,} { pp.~}                        {pages}
    +{,} { }                            {status}
    +{,} { \PrintDOI}                   {doi}
    +{,} { available at \eprint}        {eprint}
    +{}  { \parenthesize}               {language}
    +{}  { \PrintTranslation}           {translation}
    +{;} { \PrintReprint}               {reprint}
    +{.} { }                            {note}
    +{.} {}                             {transition}
}
\BibSpec{article}{%
    +{}  {\PrintAuthors}                {author}
    +{,} { \textit}                     {title}
    +{.} { }                            {part}
    +{:} { \textit}                     {subtitle}
    +{,} { \PrintContributions}         {contribution}
    +{.} { \PrintPartials}              {partial}
    +{,} { }                            {journal}
    +{}  { \textbf}                     {volume}
    +{}  { \PrintDatePV}                {date}
    +{,} { \eprintpages}                {pages}
    +{,} { }                            {status}
    +{,} { \PrintDOI}                   {doi}
    +{,} { available at \eprint}        {eprint}
    +{}  { \parenthesize}               {language}
    +{}  { \PrintTranslation}           {translation}
    +{;} { \PrintReprint}               {reprint}
    +{.} { }                            {note}
    +{.} {}                             {transition}
}
\BibSpec{book}{%
    +{}  {\PrintPrimary}                {transition}
    +{,} { \textit}                     {title}
    +{.} { }                            {part}
    +{:} { \textit}                     {subtitle}
    +{,} { \PrintEdition}               {edition}
    +{}  { \PrintEditorsB}              {editor}
    +{,} { \PrintTranslatorsC}          {translator}
    +{,} { \PrintContributions}         {contribution}
    +{,} { }                            {series}
    +{,} { \voltext}                    {volume}
    +{,} { }                            {publisher}
    +{,} { }                            {organization}
    +{,} { }                            {address}
    +{,} { pp.~}                        {pages}
    +{,} { \PrintDateB}                 {date}
    +{,} { }                            {status}
    +{}  { \parenthesize}               {language}
    +{}  { \PrintTranslation}           {translation}
    +{;} { \PrintReprint}               {reprint}
    +{.} { }                            {note}
    +{.} {}                             {transition}
}
\fi
\newcommand\go{G^{(0)}} 
\newcommand\ho{H^{(0)}}
\newcommand\lo{L^{(0)}}
\renewcommand\lg{\lambda} 
\newcommand\lh{\beta}
\renewcommand\ll{\kappa}
\def\g[#1,#2]{{}_{G}[#1,#2]} 
\def\h[#1,#2]{[#1,#2]_{H}}
\newcommand\op{\text{\normalfont op}}
\newcommand\zop{Z^{\op}}
\newcommand\cc{C_{c}}
\newcommand\pg{p_{G}}
\newcommand\ph{p_{H}}
\newcommand\mydot{\mathbin{:}}
\newcommand\Lip{\lip \scriptstyle\star}
\newcommand\Rip{\rip \scriptstyle\star}
\newcommand\tLip{\tlip \scriptstyle\star}
\newcommand\tRip{\trip \scriptstyle\star}
\newcommand\X{\mathsf{X}}
\newcommand\Y{\mathsf{Y}}
\newcommand\gux{G\backslash X}
\renewcommand\H{\mathcal{H}}
\DeclareMathOperator{\Lin}{Lin}
\newcommand\HH{\mathscr{H}}
\newcommand\hoo{\H_{00}}
\def\ipp(#1|#2){\ip({#1}|{#2})_{\pi}}
\newcommand\xind{\X\text{-}\Ind}
\newcommand\I{\mathscr{I}}
\newcommand\clsp{\overline{\operatorname{span}}}

\allowdisplaybreaks

\usepackage{color}
\definecolor{Dgreen}{cmyk}{0.93,0.33,0.92,0.25} 

\definecolor{Myblue}{RGB}{24,12,128}

\begin{document}

\author{Aidan Sims}
\address{School of Mathematics and Applied Statistics\\
University of Wollongong \\
NSW 2522\\
Australia}

\email{asims@uow.edu.au}

\author{Dana P. Williams}
\address{Department of Mathematics \\ Dartmouth College \\ Hanover, NH
03755-3551}

\email{dana.williams@Dartmouth.edu}

\subjclass[2000]{46L55}

\keywords{Groupoid, groupoid equivalence, reduced
$C^*$-algebra, equivalence theorem, Hilbert bimodule,
$C^*$-correspondence, Morita equivalence, Haar system}

\thanks{This research was supported by the Australian Research Council
  and the Edward Shapiro fund at Dartmouth College.}
\thanks{Aidan thanks Dana and the Department or Mathematics at Dartmouth
College for their warm hospitality and support.}

\date{16 February 2010; Revised 12 July 2010} 

\title[Renault's Equivalence Theorem]{\boldmath Renault's Equivalence
  Theorem for Reduced Groupoid \cs-algebras}

\begin{abstract}
We use the technology of linking groupoids to show that
equivalent groupoids have Morita equivalent reduced
$C^*$-algebras. This equivalence is compatible in a natural way
in with the Equivalence Theorem for full groupoid
$C^*$-algebras.
\end{abstract}

\maketitle

\section*{Introduction}\label{sec:introduction}

Renault's Equivalence Theorem is one of the fundamental tools in the
theory of groupoid \cs-algebras.  It states that if $G$ and $H$ are
equivalent via a $(G,H)$-equivalence $Z$, then the groupoid
\cs-algebras $\cs(G)$ and $\cs(H)$ are Morita equivalent via an \ib\
$\X$ which is a completion of $\cc(Z)$.  However, one is often
interested in the reduced \cs-algebras $\cs_{r}(G)$ and
$\cs_{r}(H)$. For example, it is the reduced \cs-algebras that play a
role in Baum-Connes theory.  Furthermore, it is the reduced algebra
--- rather than the full one --- which arises in many applications
because it, and its reduced norm, have much more concrete descriptions
than their universal counterparts. It is apparently ``well known'' to
experts that equivalent groupoids have Morita equivalent reduced
\cs-algebras.  For example, it is listed as a consequence of the main
result in \cite{tu:doc04} (see Corollary~7.9 of
\cite{tu:doc04}*{Theorem~7.8}).  It is also stated without proof
immediately following \cite{ren:amsts06}*{Theorem~3.1}.

The purpose of this paper is three fold: firstly to give a
precise statement and proof of the equivalence result for
reduced groupoid \cs-algebras; secondly to illustrate that the
equivalence result for reduced algebras is compatible with the
result for the full algebras and Rieffel induction; and
thirdly, and possibly most importantly, to highlight the role
of the linking groupoid, which is the main tool in our proofs.
The concept of the linking groupoid $L$ of an equivalence
between groupoids $G$ and $H$ is first alluded to at the end of
\cite{ren:pspm82}*{\S3} and appears in work of Kumjian --- see
in particular, \cite{kum:cjm86}. The linking groupoid was
described in general in Muhly's unpublished notes
\cite{muh:cbms}*{Remark~5.35}.  A missing ingredient up until
recently has been a Haar system for $L$. We show that if $G$
and $H$ have Haar systems, then so does $L$; we may then form
$\cs(L)$, and we show that it is isomorphic to the linking
algebra $L(\X)$ of Renault's \ib\ $\X$
(Corollary~\ref{cor-linking}).\footnote{Walther Paravicini has
also
  recently produced a Haar system for linking groupoids in his Ph.D.\
  thesis \cite{par:thesis07}*{Proposition~6.4.5}.  Parts of his work
  appear in \S1.6 of \cite{par:jkt10}, where he also proves results
  related to ours for Banach algebra completions of groupoid
  algebras. We also want to thank Paravicini for bringing the results
  in \cite{tu:doc04} to our attention.}%

Our main results imply that if $G$ and $H$ are equivalent
groupoids, then their reduced groupoid \cs-algebras
$\cs_{r}(G)$ and $\cs_{r}(H)$ are Morita equivalent via a
quotient $\X_{r}$ of $\X$ (Theorem~\ref{thm-main2}). Moreover,
we show that the Rieffel correspondence associated to $\X$
matches up the kernel $I_{\cs_{r}(G)}$ of the canonical
surjection of $\cs(G)$ onto $\cs_{r}(G)$ with the kernel
$I_{\cs{r}(H)}$ of the surjection of $\cs(H)$ onto
$\cs_{r}(H)$. Therefore for any representation $\pi$ of
$\cs(H)$ that factors through $\cs_{r}(H)$, the induced
representation $\xind \pi$ of $\cs(G)$ factors through
$\cs_{r}(G)$.

Our proof of the Equivalence Theorem for the universal algebras, like
existing ones, relies heavily on Renault's Disintegration Theorem
(\cite{ren:jot87}*{Proposition~4.2}) which is a highly nontrivial
result.  We have organized our work to illustrate that, by contrast,
the Morita equivalence for the reduced algebras can be proved without
invoking the Disintegration Theorem. Therefore there is a sense in
which the equivalence result for reduced \cs-algebras is a more
elementary result than the corresponding result for the universal
algebras.

We review the set up of the Equivalence Theorem from
\cite{mrw:jot87}*{\S2} in Section~\ref{sec:equiv-theor-univ}, and we
describe the linking groupoid and its Haar system in
Section~\ref{sec:linking-groupoid}.  In Section~\ref{sec:regul-repr}
we review some basic facts about regular representations and the
reduced groupoid \cs-algebra.  We spend a bit more time than strictly
necessary so as to clear up some ambiguities in the literature and to
state some results for future reference. In
Section~\ref{sec:main-theorem-part} we prove our equivalence theorem
for the reduced algebras, and then tie this in with the universal
constructs in Section~\ref{sec:univ-norm-link}.

We also include a short appendix to clarify the hypotheses necessary
for recently published proofs of the Disintegration Theorem and
generalizations.  In particular, we show that it is not always
necessary to assume the representations involved act on separable
spaces.

Because we want to be able to appeal both the original Equivalence
Theorem and the Disintegration Theorem, it is convenient, and at times
necessary, to require all our groupoids and spaces to be second
countable locally compact Hausdorff spaces.  As we are interested in
\cs-algebras associated to groupoids, all our groupoids are assumed to
have Haar systems.  By convention, all homomorphisms between
\cs-algebras are $*$-preserving, and all representations of
\cs-algebras are nondegenerate.

\section{Background}
\label{sec:equiv-theor-univ}

Throughout, $G$ and $H$ denote second countable, locally compact
Hausdorff groupoids with Haar systems $\sset{\lg^{u}}_{u\in\go}$ and
$\sset{\lh^{v}}_{v\in\ho}$, respectively,

In order to establish our notation, it will be useful to review the
statement and set-up of the Equivalence Theorem from
\cite{mrw:jot87}*{\S2}.  First, recall that if $G$ is a locally
compact groupoid, then we say that a locally compact space $Z$ is a
\emph{$G$-space} if there is a continuous, open map $r_{Z}:Z\to\go$
and a continuous map $(\gamma,z)\mapsto \gamma\cdot z$ from
$G*Z=\set{(\gamma,z)\in G\times Z:s_{G}(\gamma)=r_{Z}(z)}$ to $Z$ such
that $r_{X}(z)\cdot z=z$ for all $z$ and $(\gamma\eta)\cdot
z=\gamma\cdot(\eta\cdot z)$ for all $(\gamma,\eta)\in G^{(2)}$ with
$s_{G}(\eta)=r_{Z}(z)$.  (Hereafter we will often drop the subscripts
on all $r$ and $s$ maps and trust that the domain is clear from
context.)  The action is \emph{free} if $\gamma\cdot z=z$ implies
$\gamma=r(z)$ and \emph{proper} if the map $(\gamma,z)\mapsto
(\gamma\cdot z,z)$ is a proper map of $G*Z$ into $Z\times Z$.  Right
actions are dealt with similarly except that the structure map is
denoted by $s$ instead of $r$.

\begin{remark}
  Nowadays, many authors do not require the structure map $r_{Z}$ of a
  $G$-space $Z$ to be open. Since it is critical in the definition of
  an equivalence (see Definition~\ref{def-equi}) that both structure
  maps be open, we include the hypothesis here to avoid ambiguities.
  It was also part of the definition of $G$-action in
  \cite{mrw:jot87}.
\end{remark}

\begin{definition}
  \label{def-equi}
  Let $G$ and $H$ be locally compact groupoids. A
  \emph{$(G,H)$-equivalence} is a locally compact space $Z$ such that
  \begin{enumerate}
  \item $Z$ is a free and proper left $G$-space,
  \item $Z$ is a free and proper right $H$-space,
  \item\label{it:actions commute} the actions of $G$ and $H$ on $Z$
    commute,
  \item $r_{Z}$ induces a homeomorphism of $Z/H$ onto $\go$, and
  \item $s_{Z}$ induces a homeomorphism of $G\backslash Z$ onto $\ho$.
  \end{enumerate}
\end{definition}

If $Z$ is a $(G,H)$-equivalence, then there is a continuous map
$(y,z)\mapsto \g[y,z]$ of $Z*_{s}Z$ to $G$ uniquely determined
by $\g[y,z]\cdot z = y$ for all $(y,z) \in Z *_s Z$.  This map
induces a topological groupoid isomorphism of $(Z*_{s}Z)/H$
onto $G$. Similarly, there is a continuous map $(y,z)\mapsto
\h[y,z]$ satisfying $y\cdot\h[y,z] = z$ for all $(y,z) \in Z*_r
Z$, and this map induces an isomorphism of $G\backslash
(Z*_{r}Z)$ onto $H$. It is shown in \cite{mrw:jot87}*{\S2} that
if $Z$ is a $(G,H)$-equivalence, then $\cc(Z)$ is a
$\cc(G)\sme\cc(H)$-bimodule with actions and pre-inner products
given as follows: for $f \in C_c(G)$, $b \in C_c(H)$, and
$\phi, \psi \in C_c(Z)$,
\begin{align}
  \label{eq:7}
  f\cdot \phi(z)&=\int_{G}f(\gamma)\phi(\gamma^{-1}\cdot
  z)\,d\lg^{r(z)}(\gamma),
  \\
  \label{eq:8} \phi\cdot
  b(z)&=\int_{H}\phi(z\cdot\eta)b(\eta^{-1})\,d\lh^{s(z)}(\eta),\\ \label{eq:9}
  \Rip<\phi,\psi>(\eta)&=\int_{G}\overline{\phi(\gamma^{-1}\cdot
    z)}\psi(\gamma^{-1}\cdot z\cdot\eta)\,d\lg^{r(z)}(\gamma) \\
  \intertext{for any $z \in Z$ such that $s(z) = r(\eta)$, and}
  \Lip<\phi,\psi>(\gamma)&=\int_{H}\phi(\gamma\cdot
  w\cdot\eta)\overline{\psi(w\cdot
    \eta)}\,d\lh^{s(w)}(\eta)\label{eq:10}
\end{align}
for any $w \in Z$ such that $r(w) = s(\gamma)$.

The content of Renault's Equivalence Theorem
(\cite{mrw:jot87}*{Theorem~2.8}) is that $\cc(Z)$ is a pre-$\cc(G)\sme
\cc(H)$-\ib\ with respect to the universal norms on $\cc(G)$ and
$\cc(H)$, and that its completion $\X$ implements a Morita equivalence
between $\cs(G)$ and $\cs(H)$.

We define the opposite space of a $(G,H)$-equivalence $Z$ to be a
homeomorphic copy $\zop:=\{\overline{z} : z \in Z\}$ of $Z$ with the
structure of a $(H,G)$-equivalence determined by
\[
r(\overline{z}) = s(z), \quad s(\overline{z}) = r(z), \quad\eta \cdot
\overline{z} := \overline{z\cdot \eta^{-1}} \quad\text{and}\quad
\overline{z}\cdot\gamma = \overline{\gamma^{-1}\cdot z};
\]
and then $C_c(Z^{\op})$ becomes a pre-$C_c(H)\sme
C_c(G)$-imprimitivity bimodule as above. For $\psi \in C_c(Z^{\op})$,
define $\psi^* \in C_c(Z)$ by $\psi^*(z) :=
\overline{\psi(\overline{z})}$. The map $\psi \mapsto \psi^*$
determines an isomorphism from the $\cs(H)\sme\cs(G)$-\ib\ completion
of $\cc(\zop)$ to the dual module $\widetilde \X$ defined in
\cite{rw:morita}*{pp. 49--50}.

Since we will sometimes use the bimodules $\cc(Z)$ and $\cc(\zop)$ in
close proximity, we will write $\psi\mydot f$ and $b\mydot \psi$ for
the right and left actions on $\cc(\zop)$, respectively, and
$\tRip<\cdot,\cdot>$ and $\tLip<\cdot,\cdot>$ for the right and left
inner products on $\cc(\zop)$, respectively.

  We should mention that there are ``one-sided'' versions of
  the equivalence theorems in the literature.  Stadler and O'uchi
  \cite{Staouc:jot99} present a definition of a correspondence $Z$
  from $G$ to $H$ which implies $C_{c}(Z)$ can be completed to a
  $\cs_{r}(G)\sme\cs_{r}(H)$-correspondence $\mathsf{Y}$
  \cite{Staouc:jot99}*{Theorem~1.4}.  That is, $\mathsf{Y}$ is a
  right-Hilbert $\cs_{r}(H)$-module and there is a homomorphism of
  $\cs_{r}(G)$ into the adjointable operators
  $\mathcal{L}(\mathsf{Y})$ on $\Y$.  (A correspondence is also known
  as a right-Hilbert bimodule.)  A $(G,H)$-equivalence is an example
  of a Stadler-O'uchi correspondence.  The Stadler-O'uchi approach was
  generalized considerably by Tu in \cite{tu:doc04}, and Tu's work
  incorporates locally Hausdorff groupoids.  As mentioned in the
  introduction, the equivalence result for the reduced algebras should
  be a consequence of his work and the functorality of the
  constructions, although few details are given (see
  \cite{tu:doc04}*{Remark~7.17}).  In addition to the Stadler-O'uochi
  and Tu approaches, Renault has another definition of a
  correspondence $Z$ from $G$ to $H$ in
  \cite{ren:amsts06}*{Definition~2.5} which also extends the notion of
  equivalence. Nevertheless, we believe the linking groupoid approach
  developed in the next section has wider applications.  In
  particular, our results show that the equivalence theorem for the
  reduced algebras is a quotient of the result for the full crossed
  products.

\section{The linking groupoid}
\label{sec:linking-groupoid}

\begin{lemma}
  \label{lem-linking-groupoid} Suppose that $G$ and $H$ are locally
  compact Hausdorff groupoids and that $Z$ is a $(G,H)$-equivalence.
  Let $L$ be the topological disjoint union
  \begin{equation*}
    L=G\sqcup Z\sqcup Z^{\op}\sqcup H,
  \end{equation*}
and let $L^0 := G^0 \sqcup H^0 \subset L$. Define $r,s : L \to
L^0$ to be the maps inherited from the range and source maps on
$G$, $Z$, $Z^{\op}$ and $H$. Let $L^{(2)} := \{(k,l) \in L
\times L : s(k) = r(l)\}$, and let $(k,l) \mapsto kl$ be the
map from $L^{(2)}$ to $L$ which restricts to multiplication on
$G$ and $H$ and to the actions of $G$ and $H$ on $Z$ and
$Z^{\op}$, and satisfies
\begin{equation*}
z\overline{y} := {_G}[z,y]\quad\text{for $(z,y) \in Z *_s Z$}\quad\text{and}
\quad\overline{y}z := [y,z]_H\quad\text{for $(y,z) \in Z *_r Z$.}
\end{equation*}
Define $l \mapsto l^{-1}$ to be the map from $L$ to $L$ which
restricts to inversion on $G$ and $H$ and satisfies
$z^{-1} = \overline{z}$ and $\overline{z}^{\;-1} = z$ for $z
\in Z$. Under these operations, $L$ is a locally compact
Hausdorff groupoid, called the \emph{linking groupoid of $Z$}.
\end{lemma}

\begin{proof}
The inverse map is clearly an involution.  Since $\h[z,z]=s(z)$ and
$\g[z,z]=r(z)$, it is easy to see that the formulas for $r$ and $s$
are satisfied.

The continuity of the inverse map follows from the continuity of the
inverse maps on $G$ and $H$ together with the definition of the
topology on $Z^{\op}$.  The continuity of multiplication follows
from continuity of multiplication in $G$ and $H$, the
continuity of the actions of $G$ and $H$ on $Z$ and $\zop$, and the
continuity of $(y,z)\mapsto \g[y,z]$ and $(y,z)\mapsto \h[y,z]$.

The associativity of multiplication follows from routine
calculations using the associativity of the groupoid operations
and actions, and property~(\ref{it:actions commute}) of the
definition of groupoid equivalence. For example, if $x, y, z
\in Z$ with $s(x) = s(y)$ and $r(y) = r(z)$, then
\begin{align*}
(x\overline y)z
    &= \g[x,y]\cdot z = \g[x,y]\cdot (y\cdot \h[y,z]) \\
    &= (\g[x,y]\cdot y)\cdot \h[y,z] = x\cdot \h[y,z]\\
    &= x(\overline y z).\qedhere
\end{align*}
\end{proof}

Given a $(G,H)$-equivalence $Z$, the range map on $Z$ induces a
homeomorphism from the orbit space $Z/H$ to $\go$.  Thus if
$u\in\go$ and $z \in Z$ with $r(z)=u$, there is a Radon measure
$\sigma_{Z}^{u}$ on $Z$, supported on the orbit $z\cdot H$,
determined by
\begin{equation}
  \sigma_{Z}^{u}(\phi)=\int_{H}\phi(z\cdot
  \eta)\,d\lh^{s(z)}(\eta)\quad\text{for $\phi\in C_{c}(Z)$.}\label{pg:sigmas}
\end{equation}
As the notation suggests, $\sigma_{Z}^{u}$ does not depend on
the choice of $z \in r^{-1}(u)$: if $y \in Z$ with $r(y)=u$
also, then $y=z\cdot \eta'$ for some $\eta'\in H$ with
$r(\eta')=s(z)$, so left-invariance of $\beta$ gives
\begin{equation*}
  \int_{H}\phi(z\cdot\eta)\,d\lh^{s(z)}(\eta)=\int_{H} \phi(z\cdot
  \eta'\eta) \,d\lh^{s(\eta')}(\eta)=\int_{H}\phi(y\cdot
  \eta)\,d\lh^{s(y)}(\eta).
\end{equation*}
Fix $\phi\in C_{c}(Z)$. By
\cite{mrw:jot87}*{Lemma~2.9(b)}, the map $z\cdot H \mapsto
\int_H \phi(z \cdot \eta)\,d\lh^{s(z)}(\eta)$ is continuous on
$Z/H$. Since $r$ induces a homeomorphism of $Z/H$ onto $\go$,
it follows that there is a continuous function on $C_c(\go)$
given by
\begin{equation*}
  u\mapsto \int_{Z}\phi(z)\,d\sigma_{Z}^{u}(z).
\end{equation*}

By symmetry, we can also define a family of measures
$\sigma_{Z^{\op}}^{v}$ on $Z^{\op}$ with $\supp
\sigma_{Z^{\op}}^{v} = r_{\zop}^{-1}(v)$.

\begin{lemma}
  \label{lem-haar-on-L}
  For each $w\in\lo$, let $\ll^{w}$ be the Radon measure on $L$ given
  on $F\in C_{c}(L)$ by
  \begin{equation*}
    \ll^{w}(F)=
    \begin{cases}
      \lg^{w}(F\restr G)+\sigma^{w}_{Z}(F\restr Z)&\text{if $w\in \go$,
        and} \\
\sigma_{Z^{\op}}^{w}(F\restr{Z^{\op}}) + \lh^{w}(F\restr H)&\text{if
  $w\in \ho$.}
    \end{cases}
  \end{equation*}
Then $\set{\ll^{w}}_{w\in\lo}$ is a Haar system for $L$.
\end{lemma}
\begin{proof}
It is clear that $\supp \ll^{w}$ is $r^{-1}(w)=L^{w}$.
Continuity follows from continuity of $\sigma_{Z}$ and
$\sigma_{Z^{\op}}$ and of the Haar systems $\lg$ and $\lh$.  It
only remains to check left invariance.

Thus, we need to establish that for $k\in L$,
\begin{equation*}
  \int_{L}F(l)\,d\ll^{r(k)}(l) = \int_{L}F(kl)\,d\ll^{s(k)}(l).
\end{equation*}
For convenience, assume that $r(k)\in\go$. (The case where
$r(k)\in\ho$ is similar.) There are two possibilities: $k\in
G$, or $k\in Z$. First suppose $k\in G$. Then for any $z$
satisfying $r(z)=s(k)$,
\begin{align*}
  \int_{L}F(k l)\,d\ll^{s(k)}(l) &=
  \int_{G}F(k\gamma)\,d\lg^{s(k)}(\gamma) + \int_{H} F(k\cdot z \cdot
  \eta) \,d\lh^{s(z)}(\eta) \\
&= \int_{G}F(\gamma)\,d\lg^{r(k)}(\gamma) + \int_{H} F\bigl((k\cdot
z)\cdot \eta\bigr)\,d\lh^{s(k\cdot z)}(\eta) \\
&= \int_{G}F(\gamma)\,d\lg^{r(k)}(\gamma) + \int_{Z}F(w)\,d\sigma^{r(k)}(w) \\
&=\int_{L}F(l)\,d\ll^{r(k)}(l).
\end{align*}

Now suppose that $k\in Z$. Then
\begin{align*}
  \int_{L}F(kl)\,d\ll^{s(k)}(l) &= \int_{Z^{\op}}F(k\overline
  z)\,d\sigma^{s(k)}_{\zop} (\overline z) + \int_{H}F(k\cdot
  \eta)\,d\lh^{s(k)} (\eta). \\
\intertext{Since we can evaluate $\sigma^{s(k)}_{\zop}$ with
  any $\overline w$ such that $r(\overline w)=s(k)$, we may in
particular take $\overline{w} = \overline{k}$, giving}
  \int_{L}F(kl)\,d\ll^{s(k)}(l)
&= \int_{G}F\bigl(\g[k,\gamma^{-1}\cdot k]\bigr)\,d\lg^{r(k)}(\gamma)
  + \int_{Z} F(z)\,d\sigma_{Z}^{r(k)}(z). \\
\intertext{Since $\g[k,\gamma^{-1}\cdot k]=\gamma$ for all $\gamma$,
we conclude that}
  \int_{L}F(kl)\,d\ll^{s(k)}(l) &=\int_{L}F(l)\,d\ll^{r(k)}(l).\qedhere
\end{align*}
\end{proof}

We will always use the Haar system $\ll$ on $L$, so we will
henceforth write $\cs(L)$ in place of $\cs(L,\ll)$.
(Similarly, we will write $\cs(G)$ in place of $\cs(G,\lg)$ and
$\cs(H)$ in place of $\cs(H,\lh)$.)

Recall that there is a unital homomorphism $M:C_{b}(\lo)\to
M(\cs(L))$ such that for $h\in C_{b}(\lo)$ and $F\in C_{c}(L)$,
\begin{equation*}
  \bigl(M(h)F\bigr)(l)=h\bigl(r(l)\bigr)F(l)\quad\text{and}\quad
  \bigl(FM(h)\bigr) (l)=F(l)h\bigl(s(l)\bigr).
\end{equation*}
In particular, we may regard the characteristic functions $\pg$
and $\ph$ of $\go$ and $\ho$ in $C_b(\go)$ as complementary
projections in $M(\cs(L))$.

For $F\in \cc(L)$, let $F_{11}=F\restr G \in C_c(G)$,
$F_{12}=F\restr Z \in C_c(Z)$, $F_{21}=F\restr {\zop} \in
C_c(\zop)$ and $F_{22}=F\restr H \in C_c(H)$.  We view $F$ as a
matrix
\begin{equation*}
  F=
  \begin{pmatrix}
    F_{11}&F_{12}\\F_{21}&F_{22}
  \end{pmatrix}.
\end{equation*}
The involution on $C_{c}(L)$ is then given by
\begin{equation*}
  F^{*}=
  \begin{pmatrix}
    F_{11}^{*}&F_{21}^{*}\\F_{12}^{*}& F_{22}^{*}
  \end{pmatrix}
,
\end{equation*}
where $F_{11}^*$ and $F_{22}^{*}$ are the images of $F_{11}$
and $F_{22}$ under the standard involutions on $C_{c}(G)$ and
$C_{c}(H)$, while $F_{12}^{*}(\overline
z)=\overline{F_{12}(z)}$ and
$F_{21}^{*}(z)=\overline{F_{21}(\overline z)}$ for all $z \in
Z$. Straightforward computations show that the convolution
product on $\cc(L)$ is given by
\begin{align*}
  F*K&=
  \begin{pmatrix}
    F_{11}&F_{12}\\F_{21}&F_{22}
  \end{pmatrix}*
  \begin{pmatrix}
    K_{11}&K_{12}\\K_{21}&K_{22}
  \end{pmatrix} \\
&=
\begin{pmatrix}
  F_{11}*K_{11}+\tRip<F_{12}^{*},K_{21}> & F_{11}\cdot K_{12}+
  F_{12}\cdot K_{22}\\
F_{21}\mydot K_{11}+ F_{22}\mydot K_{21}& \Rip<F_{21}^{*},K_{12}>+
F_{22}*K_{22}
\end{pmatrix}\\
&=
\begin{pmatrix}
  F_{11}*K_{11}+\Lip<F_{12},K_{21}^{*}> & F_{11}\cdot K_{12}+
  F_{12}\cdot K_{22}\\
(K_{11}^{*}\cdot F_{21}^{*})^{*}+ (K_{21}^{*}\cdot F_{22}^{*})^{*}&
\Rip<F_{21}^{*},K_{12}>+
F_{22}*K_{22}
\end{pmatrix}.
\end{align*}
A routine norm calculation shows that we can identify
$C_{c}(L)$ with a dense subalgebra of the linking algebra
$L(\X)$.

\begin{lemma} \label{lem-full-projections}
The complementary projections $\pg$ and $\ph$ are full in
$M(\cs(L))$.
\end{lemma}
\begin{proof}
By symmetry, it will suffice to see that $\pg$ is full. For
$F,K \in \cc(L)$,
  \begin{equation}\label{eq:4}
    \begin{pmatrix}
      F_{11}&F_{12}\\ F_{21}&F_{22}
    \end{pmatrix}
*\pg *
\begin{pmatrix}
  K_{11}&K_{12}\\ K_{21}&K_{22}
\end{pmatrix}
=
\begin{pmatrix}
  F_{11}*K_{11}& F_{11}\cdot K_{12}\\
F_{21}\cdot K_{11} & \Rip<F_{21}^{*},K_{12}>
\end{pmatrix}
.
\end{equation}
So it suffices to see that elements of the form appearing on
the right-hand side of~\eqref{eq:4} span a dense subspace of
$C^*(L)$ in the inductive-limit topology. That elements of the
form $F_{11} * K_{11}$ span a dense subspace of $\cc(G)$ and
that elements of the form $F_{11}\cdot K_{12}$ span a dense
subspace of $\cc(Z)$ follow from the existence of an
approximate identity in $\cc(\go)$ for the left actions of
$\cc(G)$ on both itself and $\cc(Z)$ (see
\cite{mrw:jot87}*{Proposition~2.10}).  That elements of the
form $F_{21} \cdot K_{11}$ span a dense subspace of
$C_c(Z^{\op})$ follows from the corresponding property for
$\cc(H)$. That the image of $\Rip<\cdot,\cdot>$ is dense in
$C_c(H)$ follows from \cite{mrw:jot87}*{Proposition~2.10} using
standard techniques as in \cite{wil:crossed}*{p.~115} (see the
proof of \cite{mrw:jot87}*{Theorem~2.8}).\end{proof}

\begin{remark}[Our proofs of the equivalence
  theorems] \label{rem-proof-equi-thm} By
  Lemma~\ref{lem-full-projections} and
  \cite{rw:morita}*{Theorem~3.19}, to prove the
    Equivalence Theorem for the full groupoid \cs-algebras, it
  suffices to show that $\pg \cs(L) \pg \cong \cs(G)$ and similarly
  for $\cs(H)$; that is, to show that the norms on $\cs(L)$ and
  $\cs(G)$ agree on the subalgebra $\cc(G)$. Indeed, let
  $\|\cdot\|_\alpha$ be any pre-$\cs$-norm on $\cc(L)$ which is
  continuous in the inductive-limit topology. Then $\|\cdot\|_\alpha$
  is dominated by the universal norm, so the completion
  $\cs_\alpha(L)$ is a quotient of $\cs(L)$ whose multiplier algebra
  contains $C_b(L^{(0)})$.  The projections $\pg$ and $\ph$ are
  complementary full projections, and $\pg \cs_\alpha(L) \pg$ is
  isomorphic to the $\|\cdot\|_\alpha$-completion, $\cs_\alpha(G)$, of
  $\cc(G)$. A similar statement holds for $H$. Hence $\pg
  \cs_\alpha(L) \ph$, which is isomorphic to the
  $\|\cdot\|_\alpha$-completion of $\cc(Z)$, is a $\cs_\alpha(G) \sme
  \cs_\alpha(H)$-imprimitivity bimodule
  (\cite{rw:morita}*{Theorem~3.19}). So to prove the equivalence
  theorem for reduced groupoid $\cs$-algebras, it will suffice to show
  that the reduced norms on $\cs_r(L)$ and $\cs_r(G)$ agree on the
  subalgebra $\cc(G)$, and similarly for $H$.

  We will indeed prove (in Proposition~\ref{prop-universal-norm}) that
  the universal norms on $\cs(L)$ and $\cs(G)$ coincide on $\cc(G)$,
  and similarly for $H$. But our proof requires Renault's
    Disintegration Theorem \cite{muhwil:nyjm08}*{Theorem~7.8} as well
    as the basic set-up of \cite{mrw:jot87}*{Theorem~2.8}. So our
  proof of the equivalence theorem via the linking groupoid does not
  substantially simplify the original proof.

By contrast, when we show in Theorem~\ref{thm-main1} that the
reduced norms on $\cs_r(L)$ and $\cs_r(G)$ coincide on
$\cc(G)$, we require only the algebraic machinery from
\cite{mrw:jot87}*{Theorem~2.8} and the approximate identity of
\cite{mrw:jot87}*{Proposition~2.10} as required to prove
Lemma~\ref{lem-full-projections}.  In particular, our proof of
the equivalence theorem for reduced $\cs$-algebras does not
require the Disintegration Theorem.
\end{remark}

\section{Regular Representations}
\label{sec:regul-repr}

If $\mu$ is a finite Radon measure on $\go$, we can form the Radon
measure $\nu:=\mu\circ\lambda$ on $G$ given on $f\in \cc(G)$ by
\begin{equation*}
  \nu(f)=\int_{\go}\int_{G}f(\gamma)\,d\lg^{u}(\gamma)\,d\mu(u).
\end{equation*}
We write $\nu^{-1}$ for the image of $\nu$ under inversion. The
associated \emph{regular representation} $\Ind\mu$ is the
representation on $L^{2}(G,\nu^{-1})$ given by
\begin{equation*}
  (\Ind\mu)(f)\xi(\gamma)=\int_{G}f(\eta)
  \xi(\eta^{-1}\gamma)\,d\lg^{r(\gamma)}(\eta)  \quad\text{for
    $f$ and $\xi$ in $\cc(G)$.}
\end{equation*}
One can check that $\Ind\mu$ is a bounded representation of $\cs(G)$
either by appealing to the general theory of induction as in
\cite{ionwil:pams08}*{\S2}, or --- with some effort, but without
recourse to the equivalence theorem for full groupoid $C^*$-algebras
upon which \cite{ionwil:pams08}*{\S2} depends --- by verifying
directly verifying directly that $\|(\Ind\mu)(f)\|\le\|f\|_{I}$ for $f
\in \cc(G)$ and extending to the completions.

If $u\in \go$ and $\delta_{u}$ is the point mass, then the
representation $\Ind\delta_{u}$ is simply the representation of
$\cc(G)$ on $L^{2}(G_{u},\lg_{u})$ given by the convolution formula.
By definition, the \emph{reduced norm} on $\cc(G)$ is
\begin{equation*}
  \|f\|_{r}=\sup\set{\|(\Ind\delta_{u})(f)\|:u\in\go}.
\end{equation*}
So $\cs_{r}(G)$ is the quotient of $\cs(G)$ by
\begin{equation*}
  I_{\cs_{r}(G)}:=\bigcap_{u\in\go}\ker\bigl(\Ind\delta_{u}\bigr).
\end{equation*}
Alternatively, one can think of $\cs_{r}(G)$ as the completion of
$\cc(G)$ with respect to the reduced norm $\|\cdot\|_{r}$.

There is some inconsistency in the literature concerning the
definition of $\|\cdot \|_{r}$.  The definition given above coincides
with that given in \cite{anaren:amenable00}*{\S6.1} and the
unpublished notes \cite{muh:cbms}*{Definition~2.46}.  However, the
definition in Renault's original
\cite{ren:groupoid}*{Definition~II.2.8} takes the supremum over all
$\Ind\mu$.  We take a moment just to make sure everyone is talking
about the same norm (see Corollary~\ref{cor-all-regs-equal}).  Let $X$
be a second countable free and proper left $G$-space.  Then $\gux$ is
a locally compact Hausdorff space, and for each $x\in X$, the map
$\gamma\mapsto \gamma\cdot x$ is a homeomorphism of $G_{r(x)}$ onto
the orbit $G\cdot x$.  Just as for the measures $\sigma_{Z}^{u}$
defined in \eqref{pg:sigmas}, we define a Radon measure $\rho^{G\cdot
  x}$ on $X$ with support $G\cdot x$ by
\begin{equation*}
  \rho^{G\cdot x}(f)=\int_{X}f(y)\,d\rho^{G\cdot x}(y)
  :=\int_{G}f(\gamma^{-1}\cdot x)\,d\lg^{r(x)}(\gamma)\quad\text{for
    $f\in \cc(X)$.}
\end{equation*}
Our definition is independent of our choice of $x$ in its orbit by
left-invariance of the Haar system $\lg$. By
\cite{mrw:jot87}*{Proposition~2.9(b)}, the map
\begin{equation*}
  G\cdot x\mapsto \int_{X}f(y)\,d\rho^{G\cdot x}(y)
\end{equation*}
is continuous on $\gux$.  Given a finite Radon measure $\mu$ on
$\gux$, we define a Radon measure $\rho_{\mu}$ on $X$ by
\begin{equation*}
  \rho_{\mu}(f)=\int_{\gux}\int_{X} f(y)\,\rho^{G\cdot
    x}(y)\,d\mu(G\cdot x).
\end{equation*}
View $\H_{0}=\cc(X)$ as a dense subspace of $L^{2}(X,\rho_{\mu})$, and
let $\Lin(\H_{0})$ be the vector space of linear operators on
$\H_{0}$. Right multiplication under the convolution product on
$C_c(G)$ determines a homomorphism $R_{\mu}^{X}:\cc(G)\to
\Lin{C_{c}(X)}$, and some tedious computations show that $R_{\mu}^{X}$
is a homomorphism satisfying the hypotheses of Renault's
Disintegration Theorem (see
\cite{muhwil:nyjm08}*{Theorem~7.8}).\footnote{We called $R_{\mu}^{X}$
  a \emph{pre-representation} in \cite{muhwil:dm08}*{Definition~4.1}.
  See Appendix~\ref{sec:separb-hypoth-dist} for the definition and
  more details.} Hence $R_{\mu}^{X}$ is bounded and extends to a
representation of $\cs(G)$ on $L^{2}(X,\rho_{\mu})$ also denoted by
$R_{\mu}^{X}$.  Of course, the regular representations $\Ind\mu$ above
are special cases of the $R_{\mu}^{X}$ obtained by letting $X=G$.

\begin{remark}[The $\ll_{w}$]
  \label{rem-kappa-sub-w}
  We will need to use the Radon measures $\sset{\ll_{w}}_{w\in\lo}$
  on~$L$, where $\ll_{w}$ is the forward image of the measure
  $\ll^{w}$ of Lemma~3 under inversion.  It is not hard to check that
  for $F\in\cc(L)$ we have
  \begin{equation*}
    \ll_{w}(F)=
    \begin{cases}
      \lg_{w}(F\restr G)+\rho_{\zop}^{w}(F\restr{\zop}) &\text{if
        $w\in \go$, and}\\
      \rho_{Z}^{w}(F\restr Z)+\lh_{w}(F\restr H)&\text{if $w\in\ho$,}
    \end{cases}
  \end{equation*}
where we have identified $\ho$ with $G\backslash Z$, and $\go$ with
$Z/H$.
\end{remark}

\begin{example}
  \label{ex-reg-deltas} Let $\mu$ be the point mass $\delta_{G\cdot
    x_{0}}$. Then $L^{2}(X,\rho_{\mu})\cong L^{2}(G\cdot
  x_{0},\rho^{G\cdot x_{0}})$ and the homeomorphism $\gamma\mapsto
  \gamma\cdot x_{0}$ of $G_{r(x_{0})}$ onto $G\cdot x_{0}$ induces a
  unitary which intertwines $R_{\delta_{G\cdot x_{0}}}^{X}$ and $\Ind
  \delta_{r(x_{0})}$.
\end{example}

\begin{example}
  \label{ex-general} Let $X$ be any second countable free and proper
  left $G$-space, let $\mu$ be a finite Radon measure on $\gux$ and
  let $\rho^{G\cdot x}$ and $\rho_{\mu}$ be as above.  Let
  $\HH=\coprod_{G\cdot x\in\gux}L^{2}(X,\rho^{G\cdot x})$. If
  $\sset{f_{i}}$ is a countable set in $\cc(X)$ which is dense in the
  inductive-limit topology, then each $f_{i}$ defines a section of
  $\HH$ by $f_{i}(G\cdot x)(y)=f(y)$. Then
  \cite{wil:crossed}*{Proposition~F.8} implies that there is a Borel
  Hilbert bundle $(\gux)*\HH$ such that $\sset{f_{i}}$ is a
  fundamental sequence (see \cite{wil:crossed}*{Definition~F.1}) with
  the property that $L^{2}(X,\rho_{\mu})$ is isomorphic to
  $L^{2}((\gux)*\HH,\mu)$. Furthermore, the representation
  $R_{\mu}^{X}$ is equivalent to the direct integral
  \begin{equation*}
    \int^{\oplus}_{\gux}R_{\delta_{G\cdot x}}^{X}\,d\mu(G\cdot x).
  \end{equation*}
\end{example}

Part of the point of Examples \ref{ex-reg-deltas}~and \ref{ex-general}
is the following observation.
\begin{lemma}
  \label{lem-regs-all}
  If $X$ is a second countable free and proper left $G$-space and if
  $\mu$ is a finite Radon measure on $\gux$, then the representations
  $R_{\mu}^{X}$ factor through $\cs_{r}(G)$.
\end{lemma}
\begin{proof}
  Using the direct integral realization of $R_{\mu}^{X}$ in
  Example~\ref{ex-general} (and the fact that the map $r:X\to\go$ is
  surjective), we clearly have
  \begin{equation*}
    \ker R_{\mu}^{X}
    \supset \bigcap_{G\cdot x\in\gux} \ker R_{\delta_{G\cdot x}}^{X}
    = \bigcap _{x\in X}\ker \Ind \delta_{r(x)}
    = \bigcap_{u\in\go} \ker \Ind\delta_{u}
    = I_{\cs_{r}(G)}.\qedhere
  \end{equation*}
\end{proof}

Since we obtain the $\Ind\mu$ as examples of the $R_{\mu}^{X}$ (by
taking $X=G$), we obtain the following.

\begin{cor}
  \label{cor-all-regs-equal}
  Suppose $G$ is a second countable locally compact Hausdorff
  groupoid.  Then for all $f\in\cc(G)$,
  \begin{equation*}
    \|f\|_{r}=\sup\set{\|(\Ind\mu)(f)\|:\text{$\mu$ is a finite Borel
        measure on $\go$}}.
  \end{equation*}
\end{cor}
\begin{remark}
  Alternatively, we could take the supremum of the
  $\|R_{\mu}^{X}(f)\|$ ranging over all second countable free and
  proper $G$-spaces $X$, and all finite Radon measures on $\gux$.
\end{remark}

\section{The equivalence theorem for
reduced groupoid \texorpdfstring{$C^*$}{C*}-algebras}
\label{sec:main-theorem-part}

As mentioned in Remark~\ref{rem-proof-equi-thm}, now that we
have the linking groupoid together with its Haar system, the
proof that an equivalence induces a Morita equivalence of the
reduced algebras is fairly close to the surface and does not
require the full power of the equivalence result for the
universal algebras.

\begin{thm}
\label{thm-main1}
Suppose that $G$ and $H$ are second countable locally compact
Hausdorff groupoids with Haar systems as above, and suppose that
$Z$ is a $(G,H)$-equivalence. If $f \in C_c(G)$, and
\begin{equation*}
  F :=\begin{pmatrix}
    f&0\\0&0
  \end{pmatrix} \in C_c(L),
\end{equation*}
then $\|F\|_{\cs_{r}(L)} = \|f\|_{\cs_{r}(G)}$.  In particular,
the completion $\X_{r}$ of $\cc(Z)$ in the norm $\|x\| :=
\|\Rip<x,x>\|_{\cs_{r}(G)}^{1/2}$, equipped with the actions
and inner products given in \eqref{eq:7}--\,\eqref{eq:10}, is a
$\cs_{r}(G)\sme \cs_{r}(H)$-\ib{} isometrically isomorphic to
$\pg\cs_{r}(L)\ph$. Hence $\cs_{r}(G)$ and $\cs_{r}(H)$ are
Morita equivalent.
\end{thm}

\begin{remark}
In the proof of Theorem~\ref{thm-main1} we will use the notation
$\rho_{\zop}^{u}$ for  the Radon measure on $\zop$ which is the
image of $\sigma_{Z}^{u}$ on $Z$ under inversion.
Although we don't need to describe $\rho_{\zop}^{u}$ for the
proof of the theorem, for the sake of symmetry, we note that it
is the Radon measure on $\zop$ supported on $\zop_{u}$ such
that for all $\psi\in \cc(\zop)$
\begin{equation*}
  \rho_{\zop}^{u}(\psi)=\int_{H}\psi(\eta^{-1}\cdot
  \overline{z_{0}})\,d\lh^{r(\overline z)}(\eta),
\end{equation*}
for any $\overline {z_{0}}$ such that $s(\overline {z_{0}})=u$.
Thus after identifying $H\cdot \overline{z_{0}}$ with $u$,
$\rho_{\zop}^{u}$ is the measure on the free and proper left
$H$-space $\zop$ defined in Section~\ref{sec:regul-repr}.
\end{remark}

\begin{proof}
  Fix $f \in \cc(G)$ and let $F$ be the corresponding element
  of $\pg\cc(L)\pg \subset \cc(L)$. The theorem follows from
  Remark~\ref{rem-proof-equi-thm} once we
  establish that $\|F\|_{\cs_{r}(L)} = \|f\|_{\cs_{r}(G)}$.

  For $u\in\go$, we have $L_{u}=G_{u}\sqcup \zop_{u}$, where
  $\zop_{u}:= \set{\overline z\in \zop:s(\overline z)=r(z)=u}$.
  By definition, $\Ind^{L}\delta_{u}$ acts on
    $L^{2}(L_{u},\ll_{u})$.  Following Remark~\ref{rem-kappa-sub-w},
    $L^{2}(L_{u},\ll_{u})=L^{2}(G,\lg_{u})\oplus
    L^{2}(\zop,\rho_{\zop}^{u})$, and with respect to this decomposition,
    $(\Ind^{L}\delta_{u})(F)=(\Ind^{G}\delta_{u})(f)\oplus 0$.  It
  follows that
  \begin{align}
    \|F\|_{\cs_{r}(L)}&:= \max
    \Bigl\{\sup_{u\in\go}\|(\Ind^{L}\delta_{u})(F)\|,
    \sup_{v\in\ho}\|(\Ind^{L}\delta_{v}) (F)\|\Bigr\} \notag \\
    &= \max \Bigl\{\|f\|_{\cs_{r}(G)},
    \sup_{v\in\ho}\|(\Ind^{L}\delta_{v}) (F)\|\Bigr\}.\label{eq:2}
  \end{align}

  For $v\in\ho$, let $Z_{v} = \set{z\in Z : s(z) = v}$. Then
  $L_{v}=Z_{v}\sqcup H_{v}$.  Furthermore,
  $L^{2}(L_{v},\ll_{v})=L^{2}(Z,\rho^{v}_{Z})\oplus L^{2}(H,\lh_{v})$.
  Here $\rho_{Z}^{v}$ is the image of $\sigma_{\zop}^{v}$ under
  inversion.  It is the Radon measure on $Z$ with support $Z_{v}$
  given on $\phi\in\cc(Z)$ by
  \begin{equation*}
    \rho_{Z}^{v}(\phi)=\int_{G}\phi(\gamma^{-1}\cdot
    z_{0})\,d\lg^{r(z_{0})}(\gamma)
  \end{equation*}
for any $z_{0}\in Z$ such that $s(z_{0})=v$.  Thus, the
identification of $\ho$ and $G\backslash Z$ induced by the
source map on $Z$ carries $\rho_{Z}^{v}$ to the measure on the
free and proper $G$-space $Z$ defined in
Section~\ref{sec:regul-repr}. Hence
$(\Ind^{L}\delta_{v})(F)=R_{\delta_{G\cdot x_{0}}}^{Z}(f)\oplus
0$. By Example~\ref{ex-reg-deltas}, we have
$\|R_{\delta_{G\cdot x_{0}}}^{Z}(f)\|\le \|f\|_{\cs_{r}(G)}$.
It follows from \eqref{eq:2} that
$\|F\|_{\cs_{r}(L)}=\|f\|_{\cs_{r}(G)}$.
\end{proof}

\section{The universal norm and the linking algebra}
\label{sec:univ-norm-link}

\begin{prop}
  \label{prop-universal-norm}
  Suppose that $G$ and $H$ are second countable locally compact
  groupoids with Haar systems, and that $Z$ is a $(G,H)$-equivalence.
  Let $L$ be the linking groupoid.  If $f\in \cc(G)$ and
  \begin{equation*}
    F:=
    \begin{pmatrix}
      f&0\\0&0
    \end{pmatrix}
  \end{equation*}
  is the corresponding element of $\cc(L)$, then
  $\|F\|_{\cs(L)}=\|f\|_{\cs(G)}$.
\end{prop}
\begin{proof}
  Since every representation of $\cc(L)$ restricts to a representation
  of $\cc(G)$ (possibly on a subspace of the original representation),
  we certainly have $\|F\|_{\cs(L)}\le\|f\|_{\cs(G)}$.

  To obtain the reverse inequality, let $\pi$ be a faithful
  representation of $\cs(G)$ on $\H_{\pi}$.  By the universal
  properties of the tensor product, there is a sesquilinear form
  $\ipp(\cdot|\cdot)$ on the algebraic tensor product
  $\hoo:=C_{c}(L)*\pg\odot \H_{\pi}$ such that for $F$ and $K$ in
  $\cc(L)$ we have
  \begin{align*}
    \ipp(F*\pg\tensor \xi|K*\pg\tensor\zeta)&=
    \bip(\pi(\pg*K^{*}*F*\pg)\xi|\zeta) \\
    &= \bip(\pi\bigl(K_{11}^{*}*F_{11}+
    \tRip<K_{12},F_{21}>\bigr)\xi|\zeta).
  \end{align*}
  We want to see that $\ipp(\cdot|\cdot)$ is positive.  Fix
  $t=\sum_{i=1}^{n} F^{i}\tensor\xi_{i}\in \hoo$.  Since
  \cite{mrw:jot87}*{Theorem~2.8} applied to the $(H,G)$-equivalence
  $\zop$ implies that $\tRip<\cdot,\cdot>$ makes $\cc(\zop)$ into a
  pre-Hilbert $\cs(G)$-module, \cite{rw:morita}*{Lemma~2.65} implies
  that the matrix $M=\bigl(\tRip<F_{21}^{i},F_{21}^{j}>\bigr)_{ij}$ is
  positive in $M_{n}\bigl(\cs(G)\bigr)$.  Hence $M=D^{*}D$ for some $D
  \in M_n(C^*(G))$, so there are elements $d_{ij}\in \cs(G)$ such that
  \begin{equation*}
    \tRip<F_{21}^{i},F_{21}^{j}>=\sum_{i=1}^{n}d_{ki}^{*}d_{kj}.
  \end{equation*}
  Since $((F^{j})^{*})_{12}=(F^{j}_{21})^{*}$,
  \begin{align*}
    \ipp(t|t)&= \sum_{ij}\bip(\pi\bigl((F^{i}_{11})^{*}*F_{11}^{j} +
    \tRip<F_{21}^{i} , F_{21}^{j}>\bigr)\xi_{j}|\xi_{i})\\
    &= \sum_{ij} \bip(\pi(F_{11}^{j})\xi_{j}|{\pi(F_{11}^{i})\xi_{i}})
    +
    \sum_{ijk} \bip(\pi(d_{kj})\xi_{j}|{\pi(d_{ki})\xi_{i}})\\
    &= \Bigl(\sum_{i}\pi(F_{11}^{i})\xi_{i}\bigm|
    \sum_{i}\pi(F_{11}^{i})\xi_{i}\Bigr) + \sum_{k}\Bigl(
    \sum_{i}\pi(d_{ki})\xi_{i} \bigm| \sum_{i } \pi(d_{ki})\xi_{i}
    \Bigr) \ge0.
  \end{align*}
  Therefore $\ipp(\cdot|\cdot)$ is a pre-inner product on $\hoo$. Let
  $\mathcal{N}$ denote the subspace $\{\xi \in \hoo : \ipp(\xi|\xi) =
  0\}$. Then the Cauchy-Schwarz inequality (as in
  \cite{ped:analysis}*{\S3.1.1}) implies that $\ipp(\cdot|\cdot)$
  descends to a bona fide inner product on the quotient
  $\H_{0}=\hoo/\mathcal{N}$.  Furthermore, for each $F\in \cc(L)$, we
  can define a linear map $R(F):\hoo\to\hoo$ such that
  \begin{equation*}
    R(F)(K\tensor \xi):=F*K\tensor \xi.
  \end{equation*}
  Another application of the Cauchy-Schwarz inequality shows that
  $R(F)$ defines an operator on $\H_{0}$.  An easy calculation shows
  that
  \begin{equation}
    \label{eq:1}
    \ipp(R(F)t|t')=\ipp(t|{R(F^{*})}t')\quad\text{for $t,t'\in\hoo$.}
  \end{equation}
  Furthermore, since $\pi$ is continuous in the inductive-limit
  topology, it is not hard to see that
  \begin{equation}
    \label{eq:3}
    F\mapsto \ipp(R(F)t|t')
  \end{equation}
  is also continuous in the inductive-limit topology.  Since $\cc(L)$
  has an approximate unit for the inductive-limit topology,
  \begin{equation}
    \label{eq:5}
    \operatorname{span}\set{R(F)t:\text{$F\in\cc(L)$ and $t\in\hoo$}}
  \end{equation}
  is dense in $\hoo$.  Equations \eqref{eq:1},
  \eqref{eq:3}~and~\eqref{eq:5} imply that $R:\cc(L)\to\Lin(\H_{0})$
  satisfy the hypotheses of the Disintegration Theorem
  \cite{mrw:jot87}*{Theorem~2.8} as outlined in
  Appendix~\ref{sec:separb-hypoth-dist}, and therefore $R$ is a
  bounded representation of $\cs(L)$ on the completion $\H_{R}$ of
  $\H_{0}$.

  Since $\pi$ is faithful, it suffices to show that
  \begin{equation}\label{eq:6}
    \|R(F)\|_{\cs(L)}\ge \|\pi(f)\|=\|f\|_{\cs(G)}.
  \end{equation}

  Fix $\epsilon\in(0,\|f\|)$ and fix $\xi\in\H_{\pi}$ such that
  $\|\xi\|=1$ and $\|\pi(f)\xi\|^{2}> \|\pi(f)\|^{2}-\epsilon$.  Let
  $\set{k_{\alpha}}$ be an approximate identity in $\cc(G)$ for the
  inductive-limit topology, and let
  \begin{equation*}
    K_{\alpha}=
    \begin{pmatrix}
      k_{\alpha}&0\\0&0
    \end{pmatrix}
  \end{equation*}
  be the corresponding functions in $C_{c}(L)$.  Then, since $\pi$ is
  nondegenerate,
  \begin{equation*}
    \lim_{\alpha}\|K_{\alpha}\tensor\xi\|^{2} =
    \lim_{\alpha}\bip(\pi(k_{\alpha}^{*} *k_{\alpha})\xi|\xi)
    = \lim_{\alpha}\|\pi(k_{\alpha})\xi\|^{2}=1.
  \end{equation*}
  It follows that
  \begin{align*}
    \|R(F)\|^{2}&\ge
    \limsup_{\alpha}\|R(F)(K_{\alpha})\tensor\xi\|^{2} =
    \limsup_{\alpha}
    \bip(\pi\bigl(f^{*}*k_{\alpha}^{*}*k_{\alpha}*f\bigr)\xi|\xi)
    \\
    & = \lim_{\alpha}\|\pi(k_{\alpha})\pi(f)\xi\|^{2} =
    \|\pi(f)\xi\|^{2}>\|\pi(f)\|^{2}-\epsilon.
  \end{align*}
  Since $\epsilon$ is arbitrary, \eqref{eq:6} holds.  This completes
  the proof.
\end{proof}

As an immediate consequence of Proposition~\ref{prop-universal-norm}
and Remark~\ref{rem-proof-equi-thm}, we get the following.
\begin{cor}
  \label{cor-linking}
  Suppose that $G$ and $H$ are second countable locally compact
  groupoids with Haar systems, and that $Z$ is a $(G,H)$-equivalence.
  If $\X$ is the corresponding $\cs(G)\sme\cs(H)$-\ib\ and if $L$ is
  the linking groupoid, then $\cs(L)$ is isomorphic to the linking
  algebra $L(\X)$.
\end{cor}

Recall that if $\X$ is an $A\sme B$-\ib, then the Rieffel
correspondence provides a lattice isomorphism $\xind$ from the lattice
of ideals $\I(B)$ of $B$ and the lattice of ideals $\I(A)$ in $A$
\cite{rw:morita}*{Theorem~3.22}.  We can now prove the second part of
our main result.

\begin{thm}
  \label{thm-main2}
  Suppose that $G$ and $H$ are second countable locally compact
  groupoids with Haar systems, and that $Z$ is a $(G,H)$-equivalence.
  Let $\X$ be the associated $\cs(G)\sme\cs(H)$-\ib.  Then
  $\xind(I_{\cs_{r}(H)})=I_{\cs_{r}(G)}$.  Furthermore if $\X_{r}$ is
  the $\cs_{r}(G)\sme\cs_{r}(H)$-\ib\ of Theorem~\ref{thm-main1}, then
  the identity map from $\cc(Z)\subset \X$ to $\cc(Z)\subset \X_{r}$
  induces an isomorphism of the quotient \ib\ $\X/\X\cdot
  I_{\cs_{r}(H)}$ onto $\X_{r}$.
\end{thm}
\begin{proof}
  If $\phi\in\cc(Z)$, then
  \begin{equation*}
    \|\phi\|_{\X}^{2}
    = \|\Rip<\phi,\phi>\|_{\cs(H)}
    \ge \|\Rip<\phi,\phi>\|_{\cs_{r}(H)}
    = \|\phi\|_{\X_{r}}^{2}.
  \end{equation*}
  Therefore the identity map from $\cc(Z)\subset \X_{r}$
    to $\cc(Z) \subset \X$ induces a surjection of $\X$ onto
    $\X_{r}$. Let $Y$ denote the kernel of this surjection. Then $\Y$
  is a closed sub-bimodule of $\X$ such that $\X_{r}$ is isomorphic to
  $\X/\Y$ as \ib s.

  The Rieffel correspondence (in the form of
  \cite{rw:morita}*{Theorem~3.22} and \cite{rw:morita}*{Lemma~3.23})
  implies that
  \begin{equation*}
    \Y=\X\cdot I=J\cdot \X,
  \end{equation*}
  where $I$ and $J$ are ideals in $\cs(H)$ and $\cs(G)$,
    respectively,
  such that $\xind(I)=J$, and where
  \begin{equation*}
    I=\clsp\set{\Rip<x,y>:\text{$x\in\X$ and $y\in\Y$}} = \clsp\set{
      \Rip<y,y>:y\in Y}.
  \end{equation*}
  Thus $I\subset I_{\cs_{r}(H)}$.  On the other hand, if $b\in
  I_{\cs_{r}(H)}$, then for all $x$ and $y$ in $\X$, we have
  $\Rip<x,y>b=\Rip<x,y\cdot b>\in I$.  Since $\Rip<\cdot,\cdot>$ is
  full, it follows that $b\in I$.  Therefore
  $I=I_{\cs_{r}(H)}$. Similarly, we also must have $J=I_{\cs_{r}(G)}$.
  This completes the proof.
\end{proof}

\begin{cor}
  \label{cor-ind-red}
  Suppose that $G$, $H$ and $Z$ are as in Theorem~\ref{thm-main2}.  If
  $\pi$ is a representation of $\cs(H)$ that factors through
  $\cs_{r}(H)$, then $\xind \pi$ factors through $\cs_{r}(G)$.
\end{cor}
\begin{proof}
  By assumption, $I_{\cs_{r}(H)}\subset \ker\pi$.  But then by
  \cite{rw:morita}*{Proposition~3.24},
  \begin{equation*}
    I_{\cs_{r}(G)}=\xind (I_{\cs_{r}(H)}) \subset
    \xind(\ker\pi)=\ker(\xind\pi).\qedhere
  \end{equation*}
\end{proof}

\appendix

\section{Separability Hypotheses in the Disintegration Theorem}
\label{sec:separb-hypoth-dist}

Let $G$ be a second countable locally compact
Hausdorff\footnote{After
  replacing $\cc(G)$ with the vector space $\mathscr{C}(G)$ of
  functions generated by the functions in $\cc(V)$ for Hausdorff open
  sets $V\subset G$, the remarks in this appendix apply equally well
  to second countable locally compact, locally Hausdorff groupoids as
  studied in \cite{muhwil:nyjm08}.}  groupoid. A
\emph{pre-representation} of $\cc(G)$ on a dense subspace $\H_{0}$ of
a Hilbert space $\H$ is a homomorphism $L:\cc(G)\to\Lin(\H_{0})$ with
the following properties.
\begin{enumerate}
\item For $f\in\cc(G)$ and $h,k\in\H_{0}$,
  $\bip(L(f)h|k)=\bip(h|{L(f^{*})k})$.
\item For each $h,k\in\H_{0}$, $f\mapsto \bip(L(f)h|k)$ is continuous
  in the inductive-limit topology on $\cc(G)$.
\item The subspace $\operatorname{span}\set{L(f):\text{$f\in\cc(G)$
      and $h\in\H_{0}$}}$ is dense in $\H$.
\end{enumerate}
Renault's Disintegration Theorem implies that \emph{if $\H$ is
  separable}, then $L$ is the restriction of a representation $\bar L$
on $\H$ which is equivalent to the integrated form of a unitary
representation of $G$.  In particular, $L$ is bounded in the
$\|\cdot\|_{I}$-norm; indeed, $\|L(f)\|\le \|f\|_{I}$ for all $f\in
\cc(G)$.

Conversely, if $L$ is $\|\cdot\|_{I}$-bounded, $L$ extends to a
representation $\bar L$ via standard arguments.

Unfortunately, the hypothesis that $\H$ (or equivalently, $\H_{0}$)
have a countable dense subset was omitted from the statement of the
Disintegration Theorem in \cite{muhwil:nyjm08}*{Theorem~7.8} as well
as in its generalizations in \cite{muhwil:nyjm08}*{Theorem~7.12} and
\cite{muhwil:dm08}*{Theorem~4.13}.  Although separability was a
standing assumption in both \cite{muhwil:nyjm08} and
\cite{muhwil:dm08}, the omission of this hypothesis in the statements
of the Disintegration results was, well, misleading at best. (Note
that $\H$ must be separable if $\bar L$ is to be equivalent to the
integrated form of some unitary representation.  The later acts on a
direct integral of Hilbert spaces, and that theory only makes sense in
the presence of separability.)

\begin{remark}[Arbitrary $\H_{0}$]
  \label{rem-arb-H}
  Fortunately, in most applications, and in particular in the
  applications in this paper, we only want to invoke the
  Disintegration Theorem to show that $L$ is bounded and therefore
  extends to a bona fide representation of $\cs(L)$ on $\H$.  (That
  is, it is not necessary to show that $L$ is the integrated form of a
  unitary representation.) \emph{When this is the case, we do not need
    the hypothesis that $\H_{0}$ is separable.}  To see that $L$ is
  bounded, we just need to establish that for each $h_{0}\in\H_{0}$ of
  norm one, $\|L(f)h_{0}\|\le\|f\|_{I}$.  For this, it suffices to
  consider the restriction of $L$ to the cyclic subspace
  \begin{equation*}
    \hoo:=\set{L(f)h_{0}:f\in\cc(G)}.
  \end{equation*}
  Then $L$ defines a pre-representation
  $L_{0}:\cc(G)\to\Lin(\hoo)$. Since $G$ is second countable, $\cc(G)$
  has a countable dense set $\sset{f_{i}}$ in the inductive-limit
  topology, and the continuity condition of a pre-representation
  implies that $\sset{L(f_{i})h_{0}}$ is dense in $\hoo$.  Then the
  Disintegration Theorem applies to $L_{0}$, and
  \begin{equation*}
    \|L(f)h_{0}\|=\|L_{0}(f)h_{0}\|\le\|f\|_{I}.
  \end{equation*}
  Therefore $L$ is bounded on $\H_{0}$ and extends as claimed.
\end{remark}



\def\noopsort#1{}\def\cprime{$'$} \def\sp{^}

\end{document}